\documentclass[oneside]{amsart}
\usepackage{amssymb,latexsym,amsmath,amsthm,enumerate,mathrsfs}
\usepackage{fancyhdr}
\usepackage{color}
\usepackage{stmaryrd}

\newtheorem*{thm*}{Theorem}

\renewcommand{\u}[1]{\underline{#1}}

\newcommand{\beq}{\begin{equation}}
\newcommand{\eeq}{\end{equation}}

\definecolor{blue}{rgb}{.2,.6,.75}
\definecolor{green}{rgb}{.4,.7,.4}

\newtheorem*{theorem*}{Theorem}

\newtheorem*{Selberg2*}{Selberg's Announced Result}
\newtheorem{remark*}{Remark}

\title[On $r$-gaps between zeros of $\zeta(s)$]
 {On $r$-gaps between zeros of the Riemann zeta-function} 

\author[Conrey]{J. B. Conrey}
\address{American Institute of Mathematics, 600 E. Brokaw Road, San Jose, CA 95112}
\email{conrey@aimath.org}

\author[Turnage-Butterbaugh]{C. L. Turnage-Butterbaugh}
\address{Department of Mathematics, Duke University, 120 Science Drive, Durham, NC 27708}
\email{ctb@math.duke.edu}


\begin{document}
\maketitle

\begin{abstract}
Under the Riemann Hypothesis, we prove for any natural number $r$ there exist infinitely many natural numbers $n$ such that $(\gamma_{n+r}-\gamma_n)/(2\pi r /\log \gamma_n) > 1 + \Theta/\sqrt{r}$
and $(\gamma_{n+r}-\gamma_n)/(2\pi r /\log \gamma_n) < 1 - \vartheta/\sqrt{r}$ for explicit absolute positive constants $\Theta$ and $\vartheta$, where $\gamma$ denotes an ordinate of a zero of the Riemann zeta-function on the critical line. Selberg published announcements of this result several times without proof. 
\end{abstract}

\section{Introduction}
Let $\zeta(s)$ denote the Riemann zeta-function, and let $\rho = \beta + i\gamma$ denote a nontrivial zero of $\zeta(s)$. Consider the sequence of ordinates of zeros in the upper half-plane
\[
0 < \gamma_1 \le \gamma_2 \le \ldots \le \gamma_n \le \gamma_{n+1} \le \ldots.
\]
It is well known that 
\[
N(T) := \sum_{0<\gamma\le T}1 \sim \frac{T}{2\pi}\log T, 
\]
from which it follows that the average gap between consecutive zeros is $2\pi/\log \gamma_n$.  Assuming the Riemann Hypothesis, $\beta = 1/2$ and $\gamma \in \mathbb{R}$. The result of this note is a proof of the following theorem.

\begin{theorem*}\label{Selberg}
Assuming the Riemann Hypothesis, for any natural number $r$ there exist infinitely many $n$ such that 
\[
\frac{\gamma_{n+r}-\gamma_n}{2\pi r /\log \gamma_n} > 1 + \frac{\Theta}{\sqrt{r}} \qquad \text{ and } \qquad
\frac{\gamma_{n+r}-\gamma_n}{2\pi r /\log \gamma_n} < 1 - \frac{\vartheta}{\sqrt{r}}
\]
for the absolute positive constants $\Theta=0.574271$ and $\vartheta= 0.299856$. Moreover, for $r$ sufficiently large, we may take $\Theta= \vartheta = 0.9065$.
\end{theorem*}

There are discrepancies in the literature regarding the correct statement of this result, which we hope to now clarify. In \cite[p. 199]{Sel47}, Selberg announced, without proof, that there exists an absolute positive constant $\theta$ such that for all positive integers $r$
\[
\limsup_{n\to \infty}\frac{\gamma_{n+r}-\gamma_n}{2\pi r / \log \gamma_n} > 1+ \theta \qquad \text{ and } \qquad \liminf_{n\to \infty}\frac{\gamma_{n+r}-\gamma_n}{2\pi r/ \log \gamma_n} < 1- \theta . 
\]
This statement was later updated in the Acknowledgements section of \cite{Mue82}, with the $\theta$ appearing above replaced with $\theta/\sqrt{r}$. Finally, in the errata of Volume 1 of his collected papers \cite[p. 355]{SelCollected}, Selberg clarified the correct statement of his result.
\begin{Selberg2*}
There exist an absolute positive constant $\theta$ such that for all positive integers $r$
\[
\limsup_{n\to \infty}\frac{\gamma_{n+r}-\gamma_n}{2\pi r / \log \gamma_n} > 1+ \theta r^{-\alpha}\qquad \text{ and } \qquad \liminf_{n\to \infty}\frac{\gamma_{n+r}-\gamma_n}{2\pi r / \log \gamma_n} < 1- \theta  r^{-\alpha}, 
\]
where $\alpha$ may be taken as 2/3, and if one assumes the Riemann Hypothesis as 1/2.
\end{Selberg2*}
Selberg did not give an indication of a proof for either statement, however Heath-Brown in \cite[p. 246-249]{Titchmarsh} provides an unconditional proof of Selberg's result in the case $r=1$ using the work of Fujii \cite{Fujii1} concerning the mean value of $S(t)$ in short intervals. (Note that $\pi S(t)$ is the argument of $\zeta(s)$ at the point $s = 1/2 + it.$) We remark that Heath-Brown's proof for $r=1$ shows that the result holds for a positive proportion of integers $n$. 

The goal of this note is to give a proof of Selberg's conditional result for all $r\ge1$ with explicit constants. To prove our theorem, we adapt a method developed by Conrey, Ghosh, and Gonek \cite{CGG84} on gaps between consecutive nontrivial zeros of $\zeta(s)$ in the interval $[0,T]$ for $T$ large. The method is conditional on the Riemann Hypothesis. To our knowledge, our proof is the first to appear in the literature for $r>1$.

For a fixed, positive integer $r$, let
\begin{equation}\label{mulambda}
\lambda_r := \limsup_{n\to \infty} \frac{\gamma_{n+r}-\gamma_n}{2\pi  /\log \gamma_n} \qquad \text{and} \qquad \mu_r := \liminf_{n\to \infty} \frac{\gamma_{n+r}-\gamma_n}{2\pi /\log \gamma_n}.
\end{equation}
By definition $\lambda_r \ge r$ and similarly $\mu_r \le r$, however random matrix theory predicts that $\lambda_r = \infty$ and $\mu_r =0$.  Following \cite{CGG84}, we compare averages of a well-chosen polynomial of the form
\begin{equation}\label{polynomial}
A(t) := \sum_{n\le X}\frac{a^\pm(n)}{n^{it}},
\end{equation}
where $X=T^{1-\delta}$ for some small $\delta>0$. To adapt for $r$-gaps, we set
\begin{equation*}
M_1 := \int_{T}^{2T}\left|A(t)\right|^2\,dt
\end{equation*}
and
\begin{equation*}
M_{2}(c_r) := \int_{-\pi c_r/\log T}^{\pi c_r/\log T}\sum_{T\le \gamma \le 2T}\left|A(\gamma+\alpha)\right|^2\,d\alpha,
\end{equation*}
where $c_r$ is some nonzero real number. We see that $M_{2}(c_r)$ is monotonically increasing and
\[
M_{2}(\mu_r) \le rM_1 \le M_{2}(\lambda_r).
\]
Therefore, if $M_{2}(c_r)<rM_1$ for some choice of $a^+(n)$ and $c_r$ then $\lambda_r > c_r$. Similarly, if $M_2(c_r)>rM_1$ for some choice of $a^-(n)$ and $c_r$ then $\mu_r < c_r$. 

Connecting their work to a previous result of Montgomery and Odlyzko \cite{MO}, Conrey, Ghosh, and Gonek show 
\begin{equation*}
\frac{M_{2}(c_r)}{M_1} = h^\pm(c_r) + o(1),
\end{equation*}
where $h(c_r)$ is defined by
\begin{equation}\label{hfunction}
h^{\pm}(c_r):=c_r\mp\frac{\Re\left( \displaystyle\sum_{kn\le X}\frac{a^\pm(n)\overline{a^\pm(kn)}g_{c_r}(k)\Lambda(k)}{kn}\right)}{\displaystyle\sum_{n\le X}\frac{|a^\pm(n)|^2}{n}}
\end{equation}
and
\begin{equation*}
g_{c_r}(k)=\frac{2\sin\left(\pi c_r\frac{\log k}{\log{T}}\right)}{\pi\log k}
\end{equation*}
so that $|g_{c_r}(k)| \le 2c_r/\log T$. The function $h^\pm(c_r)$ was introduced by Montgomery and Odlyzko to study gaps between consecutive zeros of $\zeta(s)$. In particular, they show that if one is able to find $c_r$ such that $h^+(c_r) < r$ then $\lambda_r >  c_r$ and such that if $h^-(c_r) > r$ then $\mu_r < c_r$. 

Letting $r=1$ in \eqref{mulambda}, it follows from our theorem that $\lambda_1 >1$ and $\mu_1 <1$. Quantitative bounds on $\lambda_1$ and $\mu_1$ have been obtained using the above approach, with different choices of $a(n)$ leading to improved results. See \cite{BMN10} and subsequently \cite{FengWu} for discussions of these choices. The best current quantitative bounds concerning gaps between consecutive zeros of the Riemann zeta function (under the assumption of the Riemann Hypothesis) are $\lambda_1 > 3.18$, due to Bui and Milinovich \cite{BM17} , and $\mu_1 <0.515396$, due to Preobrazhenskii \cite{Preo}. We note that the method employed in \cite{BM17}, which is based on the work of Hall \cite{Hall} and different from the method discussed above, is unconditional if one restricts the analysis to critical zeros.

\section{Proof of the theorem for fixed $r\ge 1$}\label{proof}
For large gaps for any fixed $r\ge 1$, we choose $a^+(n) = d_{\ell}(n)$, where $d_\ell$ is multiplicative and defined on prime powers by 
\[
d_\ell(p^m)= \frac{\Gamma(m+\ell)}{\Gamma(\ell)m!}.
\]
Fix $\ell\ge 1$. (In the proof, we will ultimately set $\ell$ to be an explicit value depending on $r$.) Similarly, for small gaps for any fixed $r\ge 1$, we choose $a^-(n) = \lambda(n)d_{\ell}(n)$, where $\lambda(n)$ denotes the Liouville function.\\

We now prove the result for large gaps for any fixed $r \ge 1$. Take $a^+(n) = d_\ell(n)$ for $\ell \ge 1$ an integer to be determined later. In this case the relevant mean-value to compute is well known:
\begin{equation*}
\sum_{n \le x} \frac{d_\ell(n)^2}{n} = C_\ell(\log x)^{\ell^2}+O((\log T)^{\ell^2-1})
\end{equation*}
for fixed $\ell \ge 1$, uniformly for $x \le T$, where $C_\ell$ is a constant which will not have an effect in our application. It is shown in \cite[p.422]{CGG84} that for this choice of $a^+(n)$, the equation  $M_2(c_r) / M_1 = h^+(c_r)+o(1)$ reduces to 
\begin{equation}\label{large}
h^+(c_r) = c_r -2\ell \int_{0}^{1}\frac{\sin(\pi c_r v(1-\delta))}{\pi v}(1-v)^{\ell^2}\,dv + O(1/\log T)
\end{equation}
where $\delta>0$ is as in \eqref{polynomial} and will be taken to be sufficiently small. To detect large gaps, we must show that $h^+(c_r) < r$ for fixed $r\ge 1$. By the previous discussion, this will imply $\lambda_r > c_r$. For example, using \eqref{large} we can compute the following table of values.

\begin{table}[h]
\centering
\begin{tabular}{ |c|c|c|c| } 
 \hline
$r$	&	$\ell$	&	$c_r$ 	&	$h^+(c_r)$\\ 
\hline
1	& 	2.2		&	2.337	&	0.99965\\
2	&	2.8		&	3.708	&	1.99937\\
3	&	3.3		&	4.994	&	2.99975\\
4	&	3.7		&	6.235	&	3.99950\\
5	&	4.0		&	7.448	&	4.99978\\
 \hline
\end{tabular}
\caption{For fixed $r$, the table gives values of $\ell,c_r$ for which $h^+(c_r) < r$, implying $\lambda_r > c_r$.}
\label{figure: table example}
\end{table}

In general, to prove large gaps of the desired shape, we show that $h^+(c_r) < r$ for fixed $r\ge 1$ and $c_r = r + \Theta\sqrt{r}$ with $\Theta>0$. We estimate the integral appearing in \eqref{large} as follows. Let 
\[
\int_{0}^{1}\frac{\sin(\pi c_r(1-\delta) v)}{\pi v}(1-v)^{\ell^2}\,dv = I_1 + I_2,
\]
where
\[I_1:= \int_{0}^{1/c_r}\frac{\sin(\pi c_r(1-\delta) v)}{\pi v}(1-v)^{\ell^2}\,dv
\qquad \text{ and } \qquad
I_2:=\int_{1/c_r}^{1}\frac{\sin(\pi c_r(1-\delta) v)}{\pi v}(1-v)^{\ell^2}\,dv.
\]
For $I_1$, we first observe that the integrand is positive in the range of integration and write
\begin{equation}\label{I1}
I_1 \ge  I_{1,a} + I_{1,b},
\end{equation}
say, where
\[I_{1,a}:= \int_{0}^{1/4c_r}\frac{\sin(\pi c_r(1-\delta) v)}{\pi v}(1-v)^{\ell^2}\,dv,
\]
\[ I_{1,b}:= \int_{1/4c_r}^{1/2c_r}\frac{\sin(\pi c_r(1-\delta) v)}{\pi v}(1-v)^{\ell^2}\,dv,
\]
and we have discarded the portion of the integral from $1/(2c_r)$ to $1/c_r$. Now we estimate $I_{1,a}$ and $I_{1,b}$.
For $I_{1,a}$, we compare $\sin(\pi c_r(1-\delta) v)$ to $2\sqrt{2}c_r(1-\delta)v$ and find 
\begin{equation}\label{I1a}
I_{1,a} \ge \int_{0}^{1/4c_r}\frac{2\sqrt{2}c_r(1-\delta)v}{\pi v}(1-v)^{\ell^2}\,dv= \frac{2\sqrt{2}c_r(1-\delta)}{\pi(\ell^2+1)}\left(1- \left(1-\frac{1}{4c_r}\right)^{\ell^2+1}\right).
\end{equation}
Similarly for $I_{1,b}$, we compare $\sin(\pi c_r(1-\delta) v)$ to $(4-2\sqrt{2})c_r(1-\delta)v$ and find 
\begin{equation}\label{I1b}
I_{1,b} \ge \frac{(4-2\sqrt{2})c_r(1-\delta)}{\pi(\ell^2+1)}\left(\left(1-\frac{1}{2c_r}\right)^{\ell^2+1}- \left(1-\frac{1}{4c_r}\right)^{\ell^2+1}\right).
\end{equation}
Thus by \eqref{I1}, \eqref{I1a}, and \eqref{I1b}, we have
\[
I_1 \ge \frac{2c_r(1-\delta)}{\pi(\ell^2+1)}\left(\sqrt{2}- (2\sqrt{2}-2)\left(1-\frac{1}{4c_r}\right)^{\ell^2+1} -  (2-\sqrt{2})\left(1-\frac{1}{2c_r}\right)^{\ell^2+1}\right).
\]
Furthermore, since $\exp(-x)\ge 1-x$ for $x\ge0$, it follows that
\begin{equation}\label{I1estimate}
I_1  \ge \frac{2c_r(1-\delta)}{\pi(\ell^2+1)}\left\{\sqrt{2}- (2\sqrt{2}-2)\exp\left(\frac{-(\ell^2+1)}{4c_r} \right) -  (2-\sqrt{2})\exp\left(\frac{-(\ell^2+1)}{2c_r} \right)\right\}.
\end{equation}

We now estimate the second integral $I_2$. Since $v\ge0$, we have 
\[
|I_2| \le \frac{1}{\pi}\int_{1/c_r}^{1}\frac{(1-v)^{\ell^2}}{v}\,dv \le \frac{1}{\pi}\int_{1/c_r}^{1}\frac{\exp(-\ell^2v)}{v}\,dv.
\]
Thus, by the change of variable $u=\ell^2v$, we find
\begin{equation}\label{I2estimate}
I_2 \ge \frac{-1}{\pi}\int_{\ell^2/c_r}^{\infty}\frac{\exp(-u)}{u}\,du.
\end{equation}
Combining the estimates in \eqref{I1estimate} and \eqref{I2estimate}, we have
\begin{equation*}\label{Hestimate}
\begin{split}
h^+(c_r) \le c_r&-2\ell\biggl\{\frac{2c_r(1-\delta)}{\pi(\ell^2+1)}\biggl(\sqrt{2}- (2\sqrt{2}-2)\exp\left(\frac{-(\ell^2+1)}{4c_r} \right)\\ &-(2-\sqrt{2})\exp\left(\frac{-(\ell^2+1)}{2c_r} \right)\biggr)- \frac{1}{\pi}\int_{\ell^2/c_r}^{\infty}\frac{\exp(-u)}{u}\,du\biggr\}  +O\left(1/\log T\right).
\end{split}
\end{equation*}
In this case, $c_r = r+\Theta\sqrt{r}$ where $\Theta>0$, and thus $c_r >1$ for any $r\ge 1$. Thus, letting $\ell= \sqrt{bc_r-1}$, where $b>1$ is a real number that will be chosen later, we have
\[
\ell = \sqrt{bc_r-1} \ge \sqrt{br}\sqrt{1-\frac{1}{b}}
\]
for any $r\ge 1$. Furthermore, since we always have $c_r > 1$, for any $r\ge 1$ it follows that
\[
\frac{\ell^2}{c_r}> b-1,
\]
and thus we may again increase the length of integration in $I_2$ to write
\[
\int_{\ell^2/c_r}^{\infty}\frac{\exp(-u)}{u}\,du < \int_{b-1}^{\infty}\frac{\exp(-u)}{u}\,du.
\]
Combining these estimates, we find  
\begin{equation*}\begin{split}
h^+(c_r) < r+ \Theta\sqrt{r} &- 2\sqrt{br}\sqrt{1-\frac{1}{b}}\biggl\{\frac{2(1-\delta)}{\pi b}\biggl(\sqrt{2}-(2\sqrt{2}-2)\exp\left(\frac{-b}{4}\right)\\
&-(2-\sqrt{2})\exp\left(\frac{-b}{2}\right)\biggr)-\frac{1}{\pi}\int_{b-1}^{\infty}\frac{\exp(-u)}{u}\,du\biggr\}  +O\left(\frac{1}{\log T}\right).
\end{split}\end{equation*}
To show $h^+(c_r)< r$ and prove the theorem, we set

\begin{equation*}\begin{split}
\Theta = \max_{b}\biggl\{2\sqrt{b}\sqrt{1-\frac{1}{b}}&\biggl(\frac{2}{\pi b}\biggl(\sqrt{2}-(2\sqrt{2}-2)\exp\left(\frac{-b}{4}\right)\\
&-(2-\sqrt{2})\exp\left(\frac{-b}{2}\right)\biggr)-\frac{1}{\pi}\int_{b-1}^{\infty}\frac{\exp(-u)}{u}\,du\biggr)\biggr\}.
\end{split}\end{equation*}
The choice $b=5.0107$ yields $\Theta=0.574271$. With $\delta$ sufficiently small and $T$ sufficiently large, these choices guarantee that $h^+(c_r)<r$, as desired.

We now prove the result for small gaps for any fixed $r\ge1$. The proof for small gaps is similar to the proof for large gaps, so we indicate the necessary changes. Take $a^-(n) = \lambda(n)d_\ell(n)$ for $\ell \ge 1$ fixed. It is given in \cite[p.422]{CGG84} that this choice of $a^-(n)$, yields
\begin{equation}\label{small}
h^-(c_r) = c_r +2\ell \int_{0}^{1}\frac{\sin(\pi c_r v(1-\delta))}{\pi v}(1-v)^{\ell^2}\,dv +O\left(1/\log T\right).\\
\end{equation}
To detect small gaps, we must show that $h^-(c_r) > r$ for fixed $r\ge 1$. By the previous discussion, this will imply $\mu_r < c_r$.  For example, using \eqref{small} we can compute the following table of values.

\begin{table}[h]
\centering
\begin{tabular}{ |c|c|c|c| } 
 \hline
$r$	&	$\ell$	&	$c_r$ 	&	$h^-(c_r)$\\ 
\hline
1	& 	1.1		&	0.5172	&	1.00012\\
2	&	1.4		&	1.126	&	2.00118\\
3	&	1.9		&	1.831	&	3.00072\\
4	&	2.3		&	2.588	&	4.00099\\
5	&	2.7		&	3.375	&	5.00116 \\
 \hline
\end{tabular}
\caption{For fixed $r$, the table gives values of $\ell,c_r$ for which $h^-(c_r) > r$, implying $\mu_r < c_r$.}
\label{figure: table example}
\end{table}

In general, to prove small gaps of the desired shape, we show that $h^-(c_r) < r$ for fixed $r\ge 1$ and $c_r = r - \Theta\sqrt{r}$ with $\Theta>0$. We estimate the integral appearing in \eqref{small} as before, however for brevity we will perform the calculation without writing $I_1$ as the sum of two integrals of equal length.\footnote{To see how these choices affect the size of $\vartheta$ here and in the large gaps setting, please refer to the remark following the proof.} We find

\begin{equation*}\begin{split}
h^-(c_r) &\ge c_r+2\ell\biggl\{ \frac{2c_r(1-\delta)}{\pi(\ell^2+1)}\biggl(1-\exp\left(\frac{-(\ell^2+1)}{c_r} \right)- \frac{1}{\pi}\int_{\ell^2/c_r}^{\infty}\frac{\exp(-u)}{u}\,du\biggr\}  +O\left(\frac{1}{\log T}\right).
\end{split}\end{equation*}
Let $\ell= \sqrt{bc_r-1}$ and $c_r = r - \vartheta\sqrt{r}$, with $\vartheta>0$. In this case, we do not always have $c_r >1$. Indeed, since $\vartheta>0$, if $r=1$ then $0< c_r <1$. However, if we require that $\vartheta \le 0.5$, the estimate 
\[
\ell = \sqrt{bc_r-1} > \sqrt{br}\sqrt{\frac{1}{2}-\frac{1}{b}}
\]
holds for any $r\ge 1$. The requirement that $\vartheta \le 0.5$ also implies
\[
\frac{\ell^2}{c_r}\ge b-2
\]
for any $r\ge 1$, and we may increase the length of integration in $I_2$ to write
\[
\int_{\ell^2/c_r}^{\infty}\frac{\exp(-u)}{u}\,du \le \int_{b-2}^{\infty}\frac{\exp(-u)}{u}\,du.
\]
Thus, requiring that $\vartheta\le0.5$, we may put these estimates together to write 

\begin{equation*}\begin{split}
h^-(c_r)&> r - \vartheta\sqrt{r}+ 2\sqrt{br}\sqrt{\frac{1}{2}-\frac{1}{b}}\bigg\{\frac{2(1-\delta)}{\pi b}\biggl(1-\exp\left(-b\right) \biggr)-\frac{1}{\pi}\int_{b-2}^{\infty}\frac{\exp(-u)}{u}\,du\biggr\}  + O\left(\frac{1}{\log T}\right).\end{split}\end{equation*}
To show $h^-(c_r)> r$ and thus prove the theorem, we set

\begin{equation*}\begin{split}
\vartheta = \max_{b}\biggl\{2\sqrt{b}\sqrt{\frac{1}{2}-\frac{1}{b}}&\biggl(\frac{2}{\pi b}(1-\exp\left(-b\right)\biggr)-\frac{1}{\pi}\int_{b-2}^{\infty}\frac{\exp{-u}}{u}\,du\biggr)\biggr\}.
\end{split}\end{equation*}

The choice $b=5.17305$ yields $\vartheta=0.299856$. (We note that the condition $\vartheta <0.5$ is satisfied.) With $\delta$ sufficiently small and $T$ sufficiently large, these choices guarantee that $h^+(c_r)<r$, as desired. 

\begin{remark*}In the argument above for large gaps, if we had not divided the remaining portion of $I_1$ into two smaller integrals and instead compared $\sin(\pi c_r(1-\delta)v)$ to $2c_r(1-\delta)v$ over the interval $[0, 1/2c_r]$, we would have ultimately found that one can take $\Theta = 0.447$. Instead, by carrying out the analysis on $I_1 > I_{1,a} + I_{1,b}$ (see \eqref{I1}) and estimating $I_{1,a}$ and $I_{1,b}$ separately, we were able to provide the stronger constant $\Theta = 0.570717$. One could thus slightly improve the absolute constant $\Theta$ by breaking up $I_1$ into smaller pieces over the interval $[0, 1/2c_r]$, and estimating each piece accordingly. For example, writing  $I_1 > I_{1,a'} + I_{1,b'}+I_{1,c'}+ I_{1,d'}$ where each integral has equal length of integration over the interval $[0, 1/2c_r]$, one can obtain $\Theta = 0.593234$, and comparing $I_1$ to  the sum of sixteen such smaller integrals $\Theta = 0.599648$. Similarly, for small gaps, comparing $I_1$ to the sum of two smaller integrals of equal length over the interval $[0, 1/2c_r]$ yields $\vartheta = 0.359222$; using sixteen smaller integrals of equal length of integration over the interval $[0, 1/2c_r]$ yields $\vartheta = 0.379674$.
\end{remark*}

\section{Proof of the theorem for $r$ sufficiently large} We can improve the constants $\Theta$ and $\vartheta$ appearing in the theorem if we take $r$ to be large. In fact, we will see that in this setting, we may take $\Theta = \vartheta =  0.9065.$

We first consider large gaps for sufficiently large $r$. Starting with \eqref{large}, to detect large gaps of the desired size, we must show that $h^+(c_r) < r$ for sufficiently large $r$ and $c_r = r + \Theta\sqrt{r}$ with $\Theta>0$. Choosing $\ell = B\sqrt{r}$, we have
\[
h^+(c_r) < c_r -2B\sqrt{r} \int_{0}^{1}\frac{\sin(\pi r v(1-\delta))}{\pi v}(1-v)^{B^2r}\,dv + O(1/\log T)
\]
for sufficiently large $r$. Making the change of variable $rv = w$, the above inequality becomes
\begin{equation}\label{argument}\begin{split}
h^+(c_r) &< c_r -2B\sqrt{r} \int_{0}^{r}\frac{\sin(\pi w(1-\delta))}{\pi w}\left(1-\frac{w}{r}\right)^{B^2r}\,dw + O(1/\log T)\\
&< c_r -2B\sqrt{r} \int_{0}^{r}\frac{\sin(\pi w(1-\delta))}{\pi w}\exp\left(-B^2w\right)\,dw + O(1/\log T)\\
&=c_r -2B\sqrt{r} \int_{0}^{\infty}\frac{\sin(\pi w(1-\delta))}{\pi w}\exp\left(-B^2w\right)\,dw -2B\sqrt{r}E(r)+ O(1/\log T), \\
\end{split}\end{equation}
where
\[
E(r) =\int_{r}^{\infty}\frac{\sin(\pi w(1-\delta))}{\pi w}\exp\left(-B^2w\right)\,dw.
\]
Note that as $r\to \infty$, $\sqrt{r}E(r) \to 0$, so for sufficiently large $r$ this term is negligible. Thus we set
\begin{equation*}\begin{split}
\Theta = \max_{B}\left\{2B\int_{0}^{\infty}\frac{\sin(\pi w)}{\pi w}\exp\left(-B^2w\right)\,dw\right\} = \max_{B}\left\{\frac{2B}{\pi}\arctan\left(\frac{\pi}{B^2} \right)\right\}.  \\
\end{split}\end{equation*}
The choice $B=1.502243.$ yields $\Theta = 0.9065$. With $\delta$ sufficiently small, $T$ and $r$ sufficiently large, these choices guarantee that $h^+(c_r) < r$.

We now consider small gaps for $r$ sufficiently large. We begin with \eqref{small} and let $\ell = B\sqrt{r-\sqrt{r}}$. If we assume $\vartheta <1$, then $r-\vartheta \sqrt{r} > r - \sqrt{r}$ for all $r$, and we have
\[
h^-(c_r) > c_r +2B\sqrt{r-\sqrt{r}} \int_{0}^{1}\frac{\sin(\pi (r-\sqrt{r})v(1-\delta))}{\pi v}(1-v)^{B^2(r-\sqrt{r})}\,dv + O(1/\log T).
\]
Using the change of variable $(r-\sqrt{r})v=w$, we follow an analogous argument as in the previous subsection and ultimately set

\begin{equation*}\begin{split}
\vartheta = \max_{B}\left\{2B\int_{0}^{\infty}\frac{\sin(\pi w)}{\pi w}\exp\left(-B^2w\right)\,dw\right\}= \max_{B}\left\{\frac{2B}{\pi}\arctan\left(\frac{\pi}{B^2} \right)\right\}.  \\
\end{split}\end{equation*}
As before, the choice $B=1.502243$ yields $\vartheta = 0.9065$. With $\delta$ sufficiently small, $T$ and $r$ sufficiently large, these choices guarantee that $h^-(c_r) > r$. \\

\noindent{\bf Acknowledgements.} Turnage-Butterbaugh was supported by the National Science Foundation Grant DMS-1440140 while in residence at the Mathematical Sciences Research Institute during the Spring 2017 semester. The authors thank D.A. Goldston and M.B. Milinovich for helpful comments on an earlier version of the article. We also thank the anonymous referee for a suggestion that led to improved constants in the main theorem.

\end{document}